\documentclass[12pt]{article} 
\usepackage[utf8]{inputenc} 

\usepackage[heightrounded,a4paper]{geometry}
\usepackage{graphicx}

\usepackage{paralist}
\usepackage{subfig}
\usepackage{sidecap}
\usepackage{array}

\usepackage{amsthm}
\usepackage{amsmath}
\usepackage{amsfonts}
\usepackage{mathtools}
\usepackage{todonotes}

\usepackage{float}
\usepackage{hyperref}
\hypersetup{
	colorlinks=false,
	linkbordercolor=gray,
	citebordercolor=gray,
	urlbordercolor=gray}%
\usepackage{cleveref}

\allowdisplaybreaks

\usepackage[switch,modulo]{lineno}
\usepackage{amssymb}
\usepackage[affil-it]{authblk}

\newtheorem{thm}{Theorem}[]
\crefname{thm}{Theorem}{Theorems} 
\Crefname{thm}{Theorem}{Theorems}

\newtheorem{lemma}[thm]{Lemma}
\crefname{lemma}{Lemma}{Lemmas} 
\Crefname{lemma}{Lemma}{Lemmas}

\newtheorem{conj}[thm]{Conjecture}
\crefname{conj}{Conjecture}{Conjectures} 
\Crefname{conj}{Conjecture}{Conjectures}

\newtheorem{define}[thm]{Definition}
\crefname{define}{definition}{definitions} 
\Crefname{define}{Definition}{Definitions}

\newtheorem{claim}{Claim}
\crefname{claim}{claim}{claims} 
\Crefname{claim}{Claim}{Claims}

\crefname{obs}{observation}{observations} 
\Crefname{obs}{Observation}{Observations}

\newtheorem{remark}[thm]{Remark}
\crefname{remark}{remark}{remarks} 
\Crefname{remark}{Remark}{Remarks}

\newtheorem{prop}[thm]{Proposition}
\crefname{prop}{Proposition}{Propositions} 
\Crefname{prop}{Proposition}{Propositions}

\newtheorem{cor}[thm]{Corollary}
\crefname{cor}{Corollary}{Corollaries} 
\Crefname{cor}{Corollary}{Corollaries}

\newcommand{\eps}{\varepsilon}
\newcommand{\E}{\mathbb{E}}

\DeclareMathOperator{\SAT}{SAT}

\usepackage[ruled,vlined]{algorithm2e}

\renewcommand{\Pr}{\mathbb{P}}

\title{Polarized random $k$-SAT}
\author{Joel Larsson Danielsson and Klas Markstr{\"o}m}
\begin{document}
\maketitle
\begin{abstract}

In this paper we study a variation of the random $k$-SAT problem, called polarized random $k$-SAT, which contains both the classical random $k$-SAT model and the random version of Monotone $k$-SAT another well-known NP-complete version of SAT.    In this model there is a polarization parameter $p$, and in half of the clauses each variable occurs negated with probability $p$ and pure otherwise, while in the other half the probabilities are interchanged. For $p=1/2$ we get the classical random $k$-SAT model, and at the other extreme we have the fully polarized model where $p=0$, or 1. Here there are only two types of clauses: clauses where all  $k$ variables occur pure, and clauses where all $k$ variables occur negated. That is, for $p=0$, and $p=1$,  we  get an instance of random \emph{monotone} $k$-SAT.

We show that the threshold of satisfiability does not decrease as $p$ moves away from $\frac{1}{2}$ and thus that the satisfiability threshold for polarized random $k$-SAT with $p\neq \frac{1}{2}$ is an upper bound on the threshold for random $k$-SAT  Hence the satisfiability threshold for random monotone $k$-SAT is at least as large as for random $k$-SAT, and we conjecture that asymptotically, for a fixed $k$,  the two thresholds coincide.

\end{abstract}

%----------------------------------------------------------------------------------------------------------------------------------------------------------------------------------------------------
\section{Introduction}
\subsection{Random $k$-SAT}
During the last decades the random $k$-SAT problem has been the focus for a large amount of work by both computer scientists, mathematicians and physicists. In the classic version of this problem we have $n$ Boolean variables and $m$ random clauses, with $k$ variables each, giving us a random $k$-CNF formula $F$.
The clauses are chosen  uniformly  at random among all the $2^k\binom{n}{k}$  possible such clauses. It is known that if $\alpha=\frac{m}{n}$ is small enough then $F$ is satisfiable w.h.p. (asymptotically as $n\to \infty$) and if $\alpha$ is large enough then $F$ is unsatisfiable w.h.p.
It is also known that the property of being satisfiable has a sharp threshold~\cite{Friedgutsharpthreshold} , and it has been conjectured~\cite{chvatal-reed} that this threshold asymptotically occurs at some density $\alpha_k$. This long been known for $k=2$~\cite{CS}, but remained an open question for $k\geq 3$. 
It was recently shown that this also holds for sufficiently large $k$~\cite{DSS}, and that the threshold density predicted by heuristic methods from statistical physics~\cite{a,b} is correct.
For small $k\geq 3$ there are so far only constant separation upper and lower bounds on the possible value of $\alpha_k$, with~\cite{Up1,Up2} giving $\alpha_3 \geq 3.52$ and and~\cite{Dubois:2000} giving $\alpha_3 \leq 4.506$. The methods from~\cite{a,b} predict that $\alpha_3=4.26675\ldots$, but here recent computational work~\cite{kmpl} indicates that this value might be too large. The heuristic methods from~\cite{a,b} also predict a highly complex geometry for the set of solutions when $\alpha$  is close  to the threshold. Some these have been confirmed for both other variations of random $k$-SAT and random colouring problems, see e.g.~\cite{ckm}. The large gap between the bounds for $\alpha_3$ demonstrates that the classical probabilistic methods has so far not been as effective in the analysis of the threshold behaviour for random $k$-SAT as one would have hoped.  As pointed out in~\cite{AcMo}, part of the reason for this is that the different signs of the variables in the clauses of the formula $F$ leads to problems for the classical second moment method.  

With these problems in mind we now introduce the polarized random $k$-SAT model with polarization $p$. In this model a random formula $F$ is generated but now we chose each random    clause in the following way. We first pick a $k$-set $C$ of variables uniformly at random. Next, we flip a fair coin,  and if it gives a tail we let each variable in our clause be negated with probability  $p$, independently of each other. Otherwise we negate with probability $1-p$. For $p=1/2$ this gives the usual random $k$-SAT model. For the \emph{fully polarized} case $p=0$, or 1, we instead get a model were roughly half the clauses contain only negated variables and the other half only pure variables. In the SAT literature a formula of this type is known as a \emph{monotone} $k$-SAT formula. 

Deciding satisfiability for montone $k$-SAT  formulae is well known to be NP-complete, this follows from Schaefer's \cite{Sch78} complexity classification of SAT problems.  For $k$-SAT is has long been known \cite{MR739601} that the problem remains NP-complete even if the number of occurrences of each variable is bounded  and recently \cite{DaDo} similar restrictions of monotone $k$-SAT have been shown to remain NP-complete.  However, the random version of monotone $k$-SAT has not been analyzed in the existing literature.  Now  our random model allows us to continuously move between the usual random $k$-SAT distribution and the distribution for random monotone $k$-SAT and study how the threshold for satisfiability  varies with the polarisation parameter $p$.  

One of the appealing properties of monotone $k$-SAT  is that satisfying assignments can be given a very clean combinatorial description.  Given a satisfying assignment for a monotone $k$-SAT formula $F$,  the set of variables which satisfy the pure clauses is disjoint from the set of variables which satisfy the negative clauses. If we partition $F$ into the two types of clauses as $F=F_1 \cup F_{2}$, then a solution for $F$ corresponds to two disjoint sets of variables $T_1$ and $T_2$ such that  $T_i$ is a transversal for the hypergraph defined by the clauses in $F_i$.  We thus have a description of a satisfying assignment for $F$ in terms of hypergraph properties more closely aligned with the classical machinery of probabilistic combinatorics.

%------------------------------------------------------------

\vspace{-1em}
\subsection{The Satisfiability threshold for Polarized $k$-SAT}

\begin{figure}[H]
    \centering
    \includegraphics[width=0.9\textwidth]{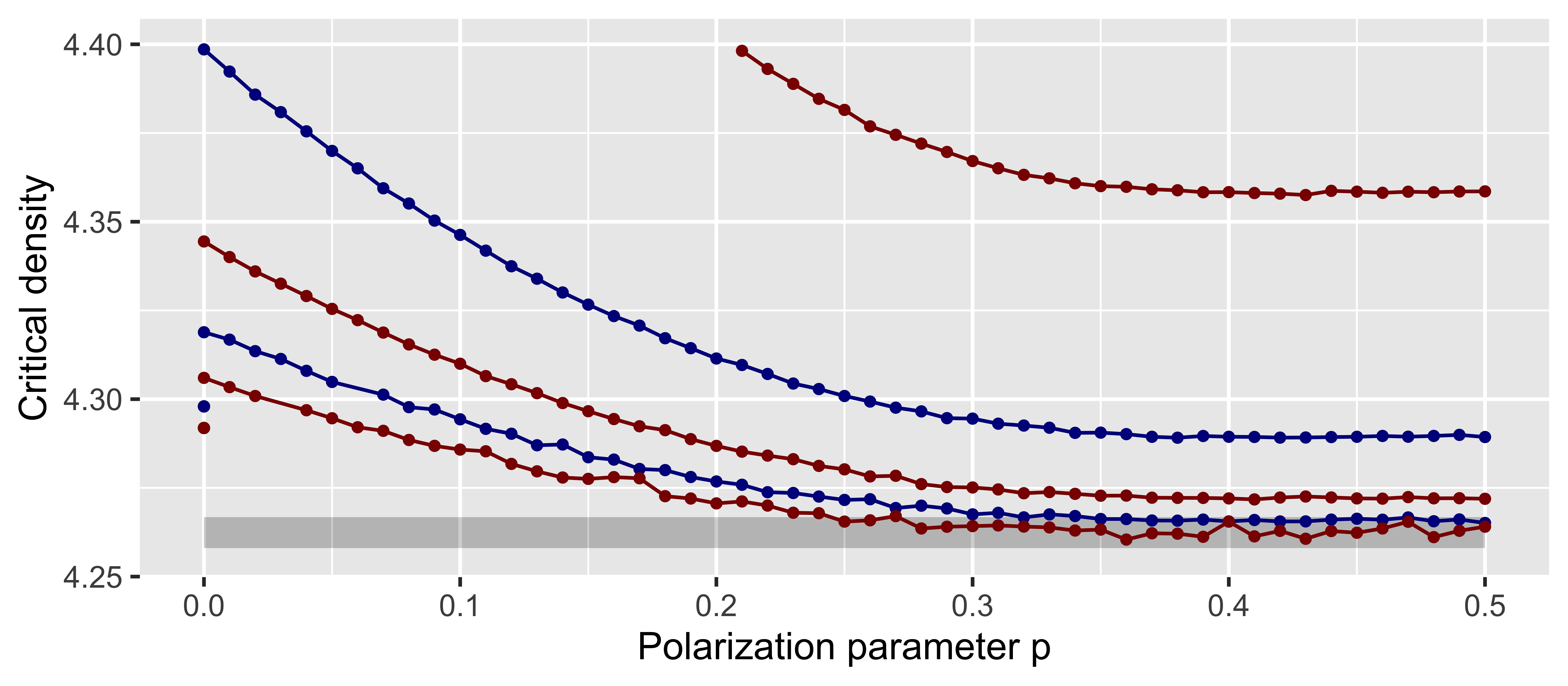}
    \caption{\textsl{Estimated critical densities of polarized $3$-SAT for several values of $n$ and $p$. From top to bottom, the curves are for $n=50$, $100$, $150$, $200$ and $250$. For $n=300$ and $350$ only simulations with $p=0$ were run; these are the two isolated points on the lower left. The shaded horizontal band is a range of predicted values for $\alpha_3$, from $4.262$ \cite{kmpl} to $4.26675$\cite{MZ2002} }}
    \label{fig:threshold}
\end{figure}

Just as for random $k$-SAT we can prove that if the  density $\alpha$ is below some constant, then a random polarized $k$-SAT formula is satisfiable with high probability. Had the clauses all been biased in the same direction 
%Had we only used a bias for negation in one direction,
(instead of having two types of clauses) we would get the \emph{biased} $k$-SAT model studied in~\cite{biasedksat}. For that model formulas become satisfiable for arbitrarily large densities, as $p$ approaches $0$. However, for polarized $k$-SAT the very general results from~\cite{flm} show that for $\alpha$ above some constant the formula is with high probability not satisfiable, independently of $p$. %Here we refine this to a $p$-dependent upper bound.
As the simulation data in \cref{fig:threshold} show, the threshold, as given by the curve for 50\% chance of satisfiability, depends on $p$ for finite $n$ and becomes flatter for larger $n$. We will study the asymptotic properties of this threshold curve.
\vspace{-1em}
The main result of this paper is that the probability of satisfiability is asymptotically non-decreasing as $p$ moves away from $\frac{1}{2}$, which is consistent with what the simulations in \cref{fig:threshold} suggest.
Our approach is to study the effect on satisfiability of adding a single clause to the formula, or switching the sign of one variable in one clause, by connecting them to the number of \emph{spine variables} of the formula using the Russo-Margulis formula and a lemma from a previous paper by the authors.

We next conjecture that for each fixed $p$ there exists an $\alpha_k(p)$ at which the model has a sharp threshold. 
The value $\alpha_k(1/2)$ is of course equal to $\alpha_k$ (if they exist), and we prove that $\alpha_k(p)$ is within a constant factor of $\alpha_k(1/2)$.
Furthermore, we conjecture that $\alpha_k(p)=\alpha_k(1/2)$, and we prove the conjecture in the special case $k=2$ by adapting a classical proof for the threshold value for random $2$-SAT.
If our second conjecture is true we thus have an alternative, and perhaps combinatorially more amenable, route to determining the threshold value for random $k$-SAT.

%------------------------------------------------------------
\vspace{-1em}
\subsection{Definitions and notation}
We will frequently use asymptotic notation to describe how functions behave as ${n\to \infty}$. $O(f), o(f), \omega(f)$ and $\Omega(f)$ will always be considered to be \emph{positive} quantities, so that we may (for instance) write $f = -O(g)$ to mean that there exists a positive constant $C$ such that $f(n) \geq -Cg(n)$ for all $n\in \mathbb{N}$. We will also use the notation $f\ll g$ to denote that $f=o(g)$. 
 
Let $\{x_i\}_{i=1}^n$ be a set of Boolean variables. We will identify TRUE (FALSE) with $+1$ ($-1$).
We say that $z$ is a \emph{literal} on the variable $x$ if $z:=x$ or $z:=\neg x$. 
A $k$-\emph{clause} is an expression of the form $z_1 \vee z_2 \vee \ldots \vee z_k$, where each $z_j$ is a literal on some variable $x_i$.
We identify $k$-clauses $C$ with the $k$-set of literals that define them, so that we can write $z_j\in C$ for the clause above. 
We say that the variable $x$ \emph{occurs} in $C$ if $C$ contains a literal on $x$. For any truth assignment $\sigma\in \{\pm 1\}^n$, we write $C(\sigma)=1$ if the clause C evaluates to TRUE when $x_i=\sigma_i$ for $i=1,2,\ldots n$, and we then say that $\sigma$ \emph{satisfies} C. Otherwise, if $C$ evaluates to FALSE,  $C(\sigma)=-1$.

A $k$-CNF (short for 'conjunctive normal form') is a Boolean formula $F$ of the form $F={C_1\wedge C_2 \wedge \ldots \wedge C_m}$, where each $C_j$ is a $k$-clause.
For any truth assignment $\sigma$, we write $F(\sigma)=1$ if $C_1(\sigma)=\ldots = C_m(\sigma)=1$ (i.e.\ all clauses are satisfied), and $F(\sigma)=-1$ otherwise (i.e.\ at least one clause is not satisfied).
If there exists a $\sigma$ such that $F(\sigma)=1$, we say that $F$ is satisfiable and write $F\in \SAT$. If no such $\sigma$ exists, we say that $F$ is unsatisfiable and write $F\notin \SAT$.

%------------------------------------------------------------
\vspace{-1em}
\subsection{The polarized $k$-SAT model}
\label{setup}
Given $0\leq p\leq 1$, we let $\Phi_m =\Phi_m(n,k,p)$ be a $p$-polarized $k$-CNF with $m$ clauses on $n$ variables. It will be convenient later to be able to separate the randomness that depends on $p$ from the randomness that doesn't. It will also be useful to couple these formulae so that $\Phi_{m-1}$ is a sub-formula of $\Phi_m$. With this in mind, we give the following more precise definition of how $\Phi_1, \Phi_2,\ldots $ are constructed.

\begin{enumerate}
\item Let $K_1,K_2,\ldots$ be a sequence of $k$-tuples $(v_1,\ldots, v_k)$ of indices in $\{1,2,\ldots, n\}$, chosen independently and uniformly at random (without replacement). 
\item Let $B_1, B_2,\ldots $ be a sequence of i.i.d.\ random variables such that \({\Pr(B_i=-1)}={\Pr(B_i=1)}=\frac{1}{2}.\)
\item For each $j=1,2,\ldots, k$, let $P_{1j},P_{2j},\ldots$ be a sequence of i.i.d.\ random variables such that $\Pr(P_{ij}=1)=p$ and $\Pr(P_{ij}=-1)=1-p$.
\item For each $k$-tuple $K_i = (v_1,\ldots, v_k)$ and each $j\leq k$, let the literal $z_j$ be defined as $x_{v_j}$ if $B_i P_{ij}=1$ and $\neg x_{v_j}$ otherwise. Then let $C_i := z_1\vee z_2\vee \ldots \vee z_k$.
\item For all $m\in \mathbb{N}$, let $\Phi_m: =C_1\wedge C_2 \wedge \ldots \wedge C_m$
\end{enumerate}

When $p=0$ or $p=1$ we say that the formula $\Phi_m$ is \emph{fully polarized}. Replacing $p$ with $1-p$ yields the same probability distribution, so we will henceforth assume (without loss of generality) that $p\leq \frac{1}{2}$. 
Note also that although the signs that the variables occur with in a clause are positively correlated\footnote{Two such variables are both pure or both negated with probability $p^2+(1-p)^2$, which is strictly greater than $\frac{1}{2}$ if $p\neq \frac{1}{2}$} for $p\neq \frac{1}{2}$, a given variable in a given clause occurs either pure or negated each with probability $\frac{1}{2}$.

Let $\alpha_k(p,n):= \min\{m/n: \Pr(\Phi_{m}\in \SAT) \leq \frac{1}{2}\}$ be the median satisfiability threshold. For $p=\frac{1}{2}$ we recover the classical random $k$-SAT problem, and we write $\alpha_k(n):= \alpha_k(\frac{1}{2},n)$.
We say that $\Phi_m$ has a \emph{sharp satisfiability threshold} if for every $\eps >0$, the formula $\Phi_m$ is satisfiable w.h.p. whenever $m<(\alpha_k(p,n)\!-\!\eps)\cdot n$ and unsatisfiable w.h.p. whenever $m>(\alpha_k(p,n)\!+\!\eps)\cdot n$.

%------------------------------------------------------------
\vspace{-1em}
\subsection{Results and conjectures}
The two main questions that we will concern ourselves with is the \emph{location} of the satisfiability threshold (as a function of $p$ or $n$), as well as its \emph{sharpness}.
Our main result is the following theorem, which lower bounds the threshold of the polarized random $k$-SAT model in terms of the classical random $k$-SAT model.
\begin{thm}
\label{mainthm}
Let $k\geq 2$ be fixed. 
\begin{enumerate}
	\item For $0\leq p \leq \frac{1}{2}$ the probability of satisfiability is asymptotically non-increasing as a function of $p$. More precisely, for any ${0\leq p \leq q\leq \frac{1}{2}}$,
\begin{equation}
\Pr(\Phi_m(n,k,p)\in \SAT)\geq \Pr(\Phi_m(n,k,q)\in \SAT)-o(1).
\label{ineq:mainthm-monotone}
\end{equation}

	\item For any ${0\leq p \leq \frac{1}{2}}$, the satisfiability threshold $\alpha_k(p,n)$ is bounded by 
\begin{equation}
\alpha_k(n)-o(1)\leq  \alpha_k(p,n) \leq  \frac{1+o(1)}{-\log_2(1-2^{-k})}.
\label{ineq:mainthm-thresholdlocation}
\end{equation}
\end{enumerate}
\end{thm}
\noindent  The first part of this theorem  tells us that asymptotically the location of the satisfiability threshold is a decreasing function of $p$, for $p\leq 1/2$. In particular for $p=0$ we get the following corollary.
\begin{cor}
	The value of satisfiability threshold for random monotone $k$-SAT  is at least as large as that for random $k$-SAT. 
\end{cor}
\noindent The second part of the theorem tells us that the classical bounds on the location of the threshold for random $k$-SAT can be generalised to polarised $k$-SAT as well.

In the special case $k=2$ of the classical random $k$-SAT model, it is well known that the threshold is sharp and $\alpha_2(n)=1+o(1)$ \cite{chvatal-reed}. Adapting this proof to the polarized model, we have shown that it too has a sharp threshold at this location.
\begin{thm}
\label{thm:2SAT}
The polarized $2$-SAT problem with polarization $p$ has a sharp satisfiability threshold and $\alpha_2(p,n)=1+o(1)$. More precisely,
\[
\Pr(\Phi_m(n,2,p)\textrm{ is satisfiable}) =
\begin{cases}
1-o(1), &m < n-\omega(n^{2/3})
\\
o(1), &m > n+\omega((n\ln n)^{9/10}),
\end{cases}
\]
In other words, $\alpha_2(p)$ exists and equals $\alpha_2 = 1$.
\end{thm}
\noindent It is known that the width of the threshold for classical random $2$-SAT is $\Theta(n^{2/3})$~\cite{2satthresholdwidth}, so it seems plausible that the lower bound in \cref{thm:2SAT} is sharp.

The location of the threshold for $k=2$ is asymptotically independent of $p$, and we conjecture that this holds for larger $k$ as well. 
\begin{conj}
\label{conj:flatthreshold}
For any $0\leq p \leq 1$ and $k\geq 3$, the threshold of the $p$-polarized random $k$-SAT problem asymptotically coincides with the threshold for the classical random $k$-SAT problem, i.e.  $\alpha_k(p,n)=\alpha_k(\frac{1}{2},n)\pm o(1)$.
\end{conj}

Friedgut~\cite{Friedgutsharpthreshold} proved a general theorem on sharp thresholds, which, as a special case, shows that the classical random $k$-SAT problem has a sharp threshold. 
Unfortunately, this theorem is not directly applicable to polarized $k$-SAT with $p\neq \frac{1}{2}$, because introducing a polarization breaks some of the symmetry of classical random $k$-SAT.\footnote{The distribution of $\Phi_m$ is not invariant under all automorpisms of the hypercube $\{-1,1\}^n$, e.g. $(x_1,x_2,\ldots,x_n) \mapsto (- x_1,x_2,\ldots,x_n)$, and Friedgut's theorem requires such symmetries.}

\begin{conj}[Generalized satisfiability conjecture]\label{conj:polarizedthreshold} Given \\ $0\leq p \leq 1$ and $k\geq 3$, there exists an $\alpha_k(p)$ such that for any $\eps>0$, \[\Pr(\Phi_m(n,k,p)\in\SAT)=\begin{cases}1-o(1), &m\leq (\alpha_k(p)-\eps)n \\ o(1), &m\geq (\alpha_k(p)+\eps)n.\end{cases}\]
\end{conj}
\noindent If \cref{conj:polarizedthreshold} is true, then $\alpha_k(p)=\lim_{n\to \infty}\alpha_k(p,n)$ exists. If \cref{conj:flatthreshold} is also true, then $\alpha_k(p)=\alpha_k(\frac{1}{2})$.

Finally, let us note that  polarised random $k$-SAT  provides a continuous interpolation between random $k$-SAT and random monotone $k$-SAT   and  our results and conjectures indicate that some properties of the satisfiability threshold do not change during this deformation. Our proof for the first part of Theorem~\ref{mainthm} depends on the structure of typical satisfying assignments for a formula, so it would be  interesting to see how existing results on the structure of the space of such assignments, of the type in e.g.  \cite{bapst_et_al:LIPIcs:2016:6645}, can be adapted to the polarised model.  Similarly it would be interesting to develop a  deterministic version of our deformation  which could unify complexity results,  like those from  \cite{MR739601}  and \cite{DaDo},  for   $k$-SAT and montone $k$-SAT.

%----------------------------------------------------------------------------------------------------------------------------------------------------------------------------------------------------
\vspace{-1em}
\section{Proof of \Cref{mainthm}}
\label{section:monotonicity}
%We will show that \cref{prop:increasingthreshold} and \cref{prop:roughbounds} together imply \Cref{mainthm}, and then we will show that \cref{prop:increasingthreshold} follows from \cref{ineq:Pprime-to-M,ineq:Pdiff-to-M} and \cref{thm:bystander} (due to Wilson~\cite{bystandervariables}). The proofs of the lemmas are postponed to \cref{subsection:proofs}.

It will be convenient to work with the parametrization $p=\frac{1}{2}-b$, where ${0\leq b\leq \frac{1}{2}}$. Let $P_{m}(b) := \Pr(\Phi_m\in\SAT)$, where $\Phi_m=\Phi_m(n,k,\frac{1}{2}-b)$.
\Cref{mainthm} will follow from the following proposition.

%The structure of the proof is as follows.
%\Cref{prop:increasingthreshold} upper bounds $P_m'$ by chaining together three inequalities: Two relating $P_m'$ (\cref{ineq:Pprime-to-M}) and $|P_m-P_{m-1}|$ (\Cref{ineq:Pdiff-to-M}) to the number of so-called \emph{spine variables}, and an upper bound on $|P_m-P_{m-1}|$ (\cref{thm:bystander}) due to Wilson~\cite{bystandervariables}.
%We use this together with the rough bounds on $P_m$ given by \cref{prop:roughbounds} to prove \cref{mainthm}.
%The proofs of the lemmas are postponed to \cref{subsection:proofs}.
\begin{prop}
\label{prop:increasingthreshold}
For any $k\geq 2$ there exists a ${c=c(k)>0}$ such that if $n>c$, $m/n\in [2^{-1},2^k]$, and $b\in[0,\frac{1}{2}]$, then
${P'_{m}(b)\leq  cn^{-\frac{k-1}{2k}}}$. 
\end{prop}
Note that $P'_m$ might well be negative, and we have proven no lower bound. If one could show that $P'_m\geq -o(1)$, this would imply \cref{conj:flatthreshold}.

\vspace{-1em}
\subsection{Spine variables and the Russo-Margulis formula}
Before we can begin with the proof of \cref{prop:increasingthreshold}, we will need some additional tools. First, the \emph{spine variables} of a random constraint satisfaction problem were introduced by Boettcher, Istrate \& Percus~\cite{spinevars} to study the computational complexity of certain algorithms. The main idea of the proof of \cref{prop:increasingthreshold} is to use spine variables in a slightly different way: A lemma from a previous paper~\cite{biasedksat} by the authors uses spine variables to characterize when a satisfiable formula can be made unsatisfiable by adding a single clause. 
We then combine this with using the \emph{Russo-Margulis formula} to study the effect on satisfiability of changing the sign of a single literal, which leads to a bound on the derivative $P_m'$.
%We then use this to relate both the derivative $P'_m$ and the difference $|P_m-P_{m-1}|$ to the number of spine variables of the formula $\Phi_{m-1}$, leading to an inequality involving $P'_m$ and $|P_m-P_{m-1}|$. Together with an upper bound on $|P_m-P_{m-1}|$ (\cref{thm:bystander}) due to Wilson~\cite{bystandervariables} this can be used to prove \cref{prop:increasingthreshold}.

Spine variables were defined in~\cite{spinevars} for both satisfiable and unsatisfiable formulae, but we will only the need the definition in the former case.
\begin{define}
Let $F$ be a satisfiable formula and $x$ a variable in it. We say that $x$ is a \emph{spine variable} in $F$ if $x$ has the same value in any assignment satisfying $F$. If such an $x$ always has value 'TRUE', we say that it is a \emph{positive spine variable} and that it is \emph{locked to TRUE}. (Similarly for negative.)
\end{define}

\begin{lemma}[From \cite{biasedksat}]
\label{spinevar:sat}
Let $F$ be a satisfiable $k$-CNF with a set $S_+\subseteq [n]$ of positive spine variables and a set $S_-\subseteq [n]$ of negative spine variables. Then $F\wedge C$ is unsatisfiable if and only if $C$ can be written as 
\[
C(x) = \Bigg(\bigvee_{i\in K_-} x_i \Bigg)\vee \Bigg(\bigvee_{i\in K_+} \neg x_i\Bigg)
\]
for some $K_\pm \subseteq S_\pm$.
\end{lemma}
In other words, $F\wedge C$ is unsatisfiable if and only if every variable in the clause $C$ is a spine variable in $F$, and its sign contradicts the value that that variable is locked to in $F$.

Next, the \emph{Russo-Margulis formula}~\cite{russosformula} is a theorem from percolation theory. It is usually stated for indicator random variables of monotone events, but we will use a slightly more general version from~\cite{grimmettpercolation}. 
\begin{thm}[Russo-Margulis, finite-dimensional case]\label{Russos-formula}
Let $I$ be a finite set, let the probability space $\mathcal{S}=\{-1,1\}^I$ be equipped with the product measure $\Pr_p$ such that if $s=\{s_i\}_{i\in I}\in \mathcal{S}$ is picked according to $\Pr_p$, then $\Pr_p(s_i=-1)=p$ for any $i \in I$. For any $s\in \mathcal{S}$ and $i\in I$, let $s^{\pm i}$ be $s$ but with the $i$:th coordinate set to $\pm 1$.
For any real random variable $X(s)$ on $\mathcal{S}$, let the \emph{pivotal} $\delta^{i}X$ be defined by
\[
\delta^{i}X(s):=
X(s^{+i})-X(s^{-i}).
\]
Then, for any $0<p<1$, 
\[
\frac{\partial }{\partial p} \E_p[X] = \sum_{i\in I} \E_p[\delta^{i}X].
\]
\end{thm}
\begin{remark}
We could have stated this result in terms of conditional expectations without defining $s^{\pm i}$ by instead writing $\E_p[\delta^i X] $ as $ \E_p[X|s_i=1]-\E_p[X|s_i=-1]$, but the definition above is more practical since it allows $X(s^{\pm i})$, $\delta^i X$, $\delta^j X$ etc.\ to all live on the same probability space.
\end{remark}

\vspace{-1em}
\subsection{Proof of \cref{prop:increasingthreshold}}
The proof consists of three inequalities, and chaining these inequalities together gives the desired result.
First, we will employ the Russo-Margulis formula with ${X=\mathbf{1}_{\Phi_m\in \SAT}}$ to upper bound $P'_{m}$.
\begin{lemma}
\label{ineq:Pprime-to-M}
Let $S$ be the (random) number of spine variables of the formula $\Phi_{m-1}$, and let ${M:=\E[S^k|\Phi_{m-1}\in \SAT]}$. Then for all large $n$ and any $b\in[0,\frac{1}{2}]$,
$P_m'(b) \leq 2k^3 P_{m-1}(b)   mn^{-k}  M^{\frac{k-1}{k}}$.
\end{lemma}
\noindent Then, we will use a similar argument to lower bound the difference $|P_m-P_{m-1}|$, again in terms of $M$. 
\begin{lemma}
\label{ineq:Pdiff-to-M}
There is a constant $c_1=c_1(k)$ such that for any $b\in[0,\frac{1}{2}]$ and $n>c_1$, $M \leq c_1 n^{k} \cdot|P_m-P_{m-1}|/P_{m-1}+c_1$.
\end{lemma}

The third inequality, which upper bounds $|P_m-P_{m-1}|$, is a result due to Wilson~\cite{bystandervariables}. They proved a lower bound on the width of the phase transition for family of random constraint satisfaction problems, including $k$-SAT. 
Crucially, this bound does not depend on the signs of variables in the random $k$-SAT formula, but only on its induced hypergraph structure. 
The following proposition is a corollary of~\cite[Theorem 1]{bystandervariables}, see their Corollary 4 for more details. We are using their result as it's stated in the last inequality of the proof of their Theorem 1.\footnote{The quantity $\eps$ in that inequality is shown to be $o(1)$ in their Lemma 3, but the proof works without modification for any $\eps\gg n^{-1}$. In particular, letting $\eps:=n^{-1/2}$ gives our stated result.}
\begin{thm}(Wilson~\cite{bystandervariables})
\label{thm:bystander}
Assume that there exist $\alpha, \alpha',\beta,N$ such that for all $n>N$ and $b\in[0,\frac{1}{2}]$, the formula $\Phi_m(n,k,\frac{1}{2}-b)$ is satisfiable with probability at least $ 1-\beta n^{-1}$ if $m<\alpha n$, and satisfiable with probability at most $ \beta n^{-1}$ if $m>\alpha' n$.

Then there exists a constant $c_2>0$ such that for all $n>c_2$, all $m_1, m_2$ such that ${\alpha < m_i/n < \alpha'}$, and all $b\in[0,\frac{1}{2}]$,
\[
|P_{m_1}(b)-P_{m_2}(b)|\leq  \frac{c_2 \cdot |m_1-m_2|}{\sqrt{n}},
\]
\end{thm}

\noindent We will also need some rough bounds on the location of the satisfiability threshold.
\begin{lemma}
\label{prop:roughbounds}
For any $k\geq 2$, there exists constants $C,c>1$ such that for any $n>C$, $p\in [0,1]$,  and $t>0$,
\begin{equation*}
\Pr(\Phi_m(n,k,p)\in\SAT) = 
    \begin{cases}
     \geq 1-30n^2/{t^3}, &m<n-t
    \\
     \leq  Cc^{-t}, &m>\frac{n}{-\log_2(1-2^{-k})}+t.
    \end{cases}
\end{equation*}
In particular, $1-o(1)\leq \alpha_k(p,n) \leq \frac{1}{-\log_2(1-2^{-k})}+o(1)$, with the $o(1)$-terms going to $0$ (as $n\to \infty$) uniformly in $p$.
\end{lemma}
\noindent The proofs of \cref{ineq:Pprime-to-M,ineq:Pdiff-to-M,prop:roughbounds} are postponed to \cref{subsection:proofs}.

\noindent \begin{proof}[Proof of \cref{prop:increasingthreshold}]
We will apply \cref{thm:bystander} with $m_1=m$, $m_2=m-1$.
By \cref{prop:roughbounds}, the assumptions of this theorem are satisfied with $\alpha=2^{-1},\alpha'=2^{k}$ and some large $\beta=\beta(k),N=N(k)$.
Combining this bound with \cref{ineq:Pdiff-to-M}, we get
\[
M \leq c_1 c_2 n^{k-\frac{1}{2}}/P_{m-1}+c_1 
\]
The first term on the right-hand side is at least  $c_1c_2
n^{k-\frac{1}{2}}\gg 1$, while the second is $O(1)$. So the first term dominates, and there must be some constant $c_3=c_3(k)$ such that
${M\leq c_3 n^{k-\frac{1}{2}}/P_{m-1}}$.
Plugging this bound on $M$ into  \Cref{ineq:Pprime-to-M} yields
\[
P'_m(b)\leq 2k^3 P_{m-1}^{1/k} \cdot \frac{m}{n^k}\cdot (c_3 n^{k-\frac{1}{2}})^\frac{k-1}{k}.
\]
Since $m<2^k n$ by assumption and $P_{m-1}^{1/k}\leq 1$, the right hand side is $O(n^{-\frac{k-1}{2k}})$, uniformly in $b$.
\end{proof}
\noindent \begin{proof}[Proof of \cref{mainthm}]
If $m/n\in [2^{-1},2^k]$, then by \cref{prop:increasingthreshold} there exists a $c>0$ such that
\[
P_m\Big(\frac{1}{2}-q\Big)-P_m\Big(\frac{1}{2}-p\Big)=\int_{\frac{1}{2}-q}^{\frac{1}{2}-p} P'_m(b)db \leq (q-p)cn^{-\frac{k-1}{2k}} = o(1)
\]
If $m/n$ is not in this interval, then by \cref{prop:roughbounds} the two terms on the left-hand side above are either both $o(1)$ or both $1-o(1)$, and hence their difference is at most $o(1)$. In either case, the inequality (\ref{ineq:mainthm-monotone}) follows.

For the lower bound in inequality (\ref{ineq:mainthm-thresholdlocation}),
let $m=m(n)$ be the largest integer such that $P_{m}(\frac{1}{2})\geq 0.6$. Then $P_{m+1}(\frac{1}{2})\leq 0.6$.
%and ${m_2}$ the smallest integer such that $P_{m_2}(\frac{1}{2})\leq 0.5$.
By \cref{thm:bystander}, ${|P_{m}(\frac{1}{2})-P_{m+1}(\frac{1}{2})|=o(1)}$, so ${P_m(\frac{1}{2})=0.6+o(1)}$.
%(Such an $m$ has to exist, since by \cref{thm:bystander} ${|P_m-P_{m-1}|=o(1)}$.)
Since classical $k$-SAT has a sharp satisfiability threshold property (Friedgut~\cite{Friedgutsharpthreshold}), and $P_m(\frac{1}{2})$ is bounded away from both $0$ and $1$ (i.e. $m$ is in the critical window), ${m/n=\alpha_k(n)-o(1)}$. By the previous inequality,
$P_m(\frac{1}{2}-p)\geq P_m(\frac{1}{2})-o(1)$, which is at least $0.6-o(1)> 0.5$ by assumption. Hence $\alpha_k(p,n)\geq m/n=\alpha_k(n)-o(1)$. The upper bound in inequality  (\ref{ineq:mainthm-thresholdlocation}) follows directly from \cref{prop:roughbounds}.
\end{proof}

\vspace{-1em}
\subsection{Proofs of lemmas}
\label{subsection:proofs}
The following simple inequality will be useful to us several times.
\begin{claim}
\label{claim:with-or-without-replacement}
If we pick $k$ elements from a set of $s$ elements uniformly at random with replacement, the probability that we pick $k$ distinct elements is at least ${(1-k/s)^k\geq 1-k^2/s}$. 
\end{claim}
\vspace{-1em}
\noindent \begin{proof}
The probability that we pick distinct elements is $\prod_{i=0}^{k-1}(1-i/s)\geq (1-k/s)^k$, and $(1-k/s)^k\geq 1-k^2/s$ follows from the convexity of $x\mapsto (1-x)^k$ for $x\leq 1$.
\end{proof}

\noindent \begin{proof}[Proof of \cref{ineq:Pprime-to-M}]
To apply the Russo-Margulis formula, we must figure out what the correct pivotals are. 
Recall that we constructed the random $k$-CNF $\Phi_m$ from the random variables $K_i$, $B_i$ and $P_{ij}$ ($1\leq i\leq m, 1\leq j\leq k$), where $K_i$ is the ordered list of variables occuring in the $i$:th clause $C_i$, and $B_i \cdot P_{ij}$ is the sign of the $j$:th variable in $C_i$.
Let $H:=((K_1,B_1),(K_2,B_2),\ldots ,(K_m,B_m))$ denote the (signed, ordered) hypergraph structure of the formula $\Phi_m$. Note that $H$ does not depend on $p$. 
We will condition on $H$, and study the effect on $X$ of switching the value of one of the random variables $P_{ij}$. This only affects the clause $C_i$, and leaves the rest of the formula $\Phi_m$ unaffected. 

Let $C_i^{+j}$ and $C_i^{-j}$ be the two clauses obtained from $C_i$ by letting the $j$\textsuperscript{th} variable in $C_i$ occur with sign $B_i$ and $-B_i$ respectively. (Instead of having signed given by $B_iP_{ij}$, as it has in $C_i$.) Furthermore, let the formula $F_{i}$ be defined by $F_{i}:=\Phi_m-C_i=\bigwedge_{\ell\neq i} C_\ell$. Then (still conditional on $H$) the pivotal is given by
\begin{align*}
\delta^{ij}X = &\mathbf{1}_{F_{i}\wedge C_i^{+j} \in \SAT} - \mathbf{1}_{F_{i}\wedge C_i^{-j} \in \SAT}
\\
= &\mathbf{1}_{F_{i}\wedge C_i^{-j} \notin \SAT} - \mathbf{1}_{F_{i}\wedge C_i^{+j} \notin \SAT}.
\end{align*}
Applying the Russo-Margulis formula to the product space $\{\pm 1\}^I$ with index set $I:=\{1,\ldots, m\}\times \{1,\ldots , k\}$ and with product measure $\Pr_p(\bullet) := \Pr(\bullet|H)$ leads to
\[
\frac{d}{dp}\Pr(\Phi_m\in \SAT|H) =\sum_{i=1}^m \sum_{j=1}^k  \E[ \delta^{ij}X|H]
\]
Taking expectations over $H$ of both sides gives us 
\[
\E\big[\frac{d}{dp}\Pr(\Phi_m\in \SAT|H)  \big]=\sum_{i=1}^m \sum_{j=1}^k  \E[ \delta^{ij}X], 
\]
and since there are only finitely many possible (signed, ordered) hypergraphs $H$ on $n$ vertices, we can exchange the order of expectation and differentiation in the left hand side. In the sum on the right hand side all terms are equal, and thus
\[
-\frac{dP_m}{db}=\frac{dP_m}{dp} =mk\E[ \delta^{m1}X],
\]
since $p=\frac{1}{2}-b$.
Let $S$ be the number of spine variables of $\Phi_{m-1}$ if this formula is satisfiable (and not defined otherwise). We will now lower bound $\E[\delta^{m1}X]$ in terms of $M:={\E[S^k|\Phi_{m-1}\in \SAT]}$.
Let $F:=F_{m}=\Phi_{m-1}$, and note that the pivotal $\delta^{m1}X$ is non-zero iff exactly one of $F\wedge C_m^{+1}$ and $F\wedge C_m^{-1}$ is satisfiable. 
That only happens if (i) $F$ is satisfiable and (ii) every variable in $C_m$ is one of the $S$ spine variables of $F$. Since $F$ is satisfiable with probability $P_{m-1}$,
\[
\E[ \delta^{m1}X] = P_{m-1}\E[ \delta^{m1}X|F\in \SAT]
\]
Let $S_+$ ($S_-$) be the set of positive (negative) spine variables (so that $|S_+|+|S_-|=S$), and let $\Gamma$ be the event $\{(i), (ii), |S_+|=s_+, |S_-|=s_-\}$ for some $s_\pm$.
Conditional on $\Gamma$, the variables of $C_m$ is a $k$-tuple $K_m$ of variables, chosen uniformly at random from the $s:=s_++s_-$ spine variables \emph{without} replacement. 
If $s$ is large compared to $k$, this $k$-tuple will typically be the same as a sample taken \emph{with} replacements.

We therefore let $K=(v_1,\ldots, v_k)$ be a random $k$-tuple chosen uniformly at random from $\{1, \ldots ,n\}$ \emph{with} replacement, and couple it to $K_m$ such that every element in $K$ is also in $K_m$.\footnote{One way to do this: let $u_1,u_2,\ldots$ be a sequence of integers picked uniformly at random from $[n]$ with replacement, and let $K:=(u_1,\ldots u_k)$, but let $K_m$ be the $k$ first \emph{distinct} $u_i$'s.} Then $K=K_m$ whenever $K$ consists of $k$ distinct elements, which by \cref{claim:with-or-without-replacement} happens with probability at least $1-k^2/s$.

For the sake of convenience, we will for any clause $C=z_1\vee z_2 \vee \ldots \vee z_k$ define the clause $-C$ as $-C:=\neg z_1\vee  \neg z_2 \vee \ldots \vee \neg z_k$.

Let $z_j$ be the literal $x_{v_j}$ if $P_{ij}=1$ and $\neg x_{v_j}$ otherwise, and let
the clauses $D_\pm$ be defined as ${D_+:= x_{v_1}\wedge z_2\wedge \ldots z_k}$ and ${D_-:= \neg x_{v_1}\wedge z_2\wedge \ldots z_k}$.
Let $C_+:=D_+$ if $B_m=1$, and $C_+:=-D_+$ otherwise, and analogously for $C_-$.
Then $C_+=C_m^{+1}$ and $C_-=C_m^{-1}$ with probability at least $1-k^2/s$, independently of $F$, so we can approximate $\E[\delta^{m1}X|\Gamma]$ with $\E[\mathbf{1}_{F\wedge C_- \notin \SAT} - \mathbf{1}_{F\wedge C_+ \notin \SAT}|\Gamma]$.
We will first show that the latter expected value is $\geq 0$, and then bound the error in the approximation.
\begin{claim}
\label{claim:positive-expected-pivotal}
Let $\tilde X:=\mathbf{1}_{F\wedge C_- \notin \SAT} - \mathbf{1}_{F\wedge C_+ \notin \SAT}$. Then $\E[\tilde X|\Gamma] \geq 0$.
\end{claim}
\noindent \begin{proof}[Proof of claim]
Let $a$ be the number such that $\frac{1}{2}+ a = \frac{s_+}{s}$, i.e. the fraction of positive spine variables.
Each literal $z_2,\ldots z_k$ is pure with probability $\frac{1}{2}+b$ and negated otherwise, while (condititional on $\Gamma$) each corresponding variable is a positive spine variable with probability $\frac{1}{2}+a$.
Hence the conditional probability that all of them disagree with the spine is
\[
\left(\Big(\frac{1}{2}-a\Big)\Big(\frac{1}{2}+b\Big)+\Big(\frac{1}{2}+a\Big)\Big(\frac{1}{2}-b\Big)\right)^{k-1}\!\!
=
\Big(\frac{1}{2}-2ab\Big)^{k-1}.
\]
The first literal of $D^-$ is $x_{v_1}$, and thus disagree with the spine with probability $\frac{1}{2}-a$. So, conditional on $\Gamma$, the probability that $F \wedge D_+$ is unsatisfiable is ${\big(\frac{1}{2}-a\big)\cdot \big(\frac{1}{2}+2ab\big)^{k-1}}$. Similarly, the first literal of $D_-$ is $\neg x_{v_1}$, which disagrees with the spine with probability $\frac{1}{2}+a$, so $F \wedge D_-$ is unsatisfiable with probability ${\big(\frac{1}{2}+a\big)}\cdot{ \big(\frac{1}{2}+2ab\big)^{k-1}}$.
Together, these observations give us that  
\begin{align}
&\quad\quad\Pr(F\wedge D_-\notin \SAT|\Gamma)-\Pr(F\wedge D_+\notin \SAT|\Gamma) \nonumber
\\
&=
\left(\Big(\frac{1}{2}\!+\!a\Big)-\Big(\frac{1}{2}\!-\!a\Big)\right)\cdot \Big(\frac{1}{2}+2ab\Big)^{k-1}= 2a\cdot \Big(\frac{1}{2}+2ab\Big)^{k-1}. \label{probdiff:Dplus}
\end{align}
If we consider $-D_+$ and $-D_-$, an analogous argument gives that
\begin{align}
\label{probdiff:Dminus}
\Pr(F\wedge -D_-\notin \SAT|\Gamma)-\Pr(F\wedge -D_+\notin \SAT|\Gamma) &= -2a\cdot \Big(\frac{1}{2}-2ab\Big)^{k-1}.
\end{align}
Since $C_\pm$ is equally likely to be $D_\pm$ as $-D_\pm$, averaging gives
\begin{align*}
\E[\tilde X|\Gamma] &= \Pr(F\wedge C_-\notin \SAT|\Gamma)-\Pr(F\wedge C_+\notin \SAT|\Gamma)  =\frac{1}{2}\big[(\ref{probdiff:Dplus})+(\ref{probdiff:Dminus})\big]
\\
&=a\cdot \underbrace{\bigg(\Big(\frac{1}{2}+2ab\Big)^{k-1}-\Big(\frac{1}{2}-2ab\Big)^{k-1}\bigg)}_{=:f(a,b)}.
\end{align*}
If $a\!=\!0$ or $b\!=\!0$, then $f(a,b)=0$ and hence $\E[\tilde X|\Gamma]=0$. Otherwise, if $b>0$ and $a\neq 0$, then $a$ and $f(a,b)$ have the same sign, and their product is positive.
\end{proof}

But as we noted earlier, this is not quite the expected value of the pivotal $\delta^{m1}X$.
The probability that the $k$-tuple $K$ (drawn with replacement) equals $K_m$ (drawn without replacement) is at least $1-k^2/s$ by \cref{claim:with-or-without-replacement}. 
But then $\delta^{m1}X=\tilde X$ with probability at least $1-k^2/s$ (independently of $\Gamma$), and since $|\delta^{m1}X-\tilde X|\leq 2$, 
\[
\E\big[ \delta^{m1}X\big| \Gamma\big]
\,\geq\,
\E[\tilde X|\Gamma]
-2\Pr\big( \delta^{m1}X\neq \tilde X\big|\Gamma\big)
\,\geq\,
-\frac{2k^2}{s}.
\]
Recall that $\Gamma = \{F\in \SAT, \textrm{all variables in $C$ are spines}, |S_+|=s_+, |S_-|=s_-\}$. The probability that all variables in $C$ are spine variables in $F$ is ${\binom{s}{k}/\binom{n}{k}\leq (s/n)^k}$, and if some variable is not, then $\delta^{m1}X=0$. Thus
\[
\E\Big[ \delta^{m1}X\Big| F\in \SAT,|S_+|\!=\!s_+,|S_-|\!=\!s_-\Big] = \frac{\binom{s}{k}}{\binom{n}{k}}\cdot \E\Big[ \delta^{m1}X\Big| \Gamma\Big] \geq -\left(\frac{s}{n}\right)^k \cdot  \frac{2k^2}{s}
\]
The conditional expectation $\E\big[ \delta^{m1}X\big| F\in \SAT, S=s\big]$ is a weighted average, over all $s_\pm$  with $s_++s_-=s$, of the LHS above. Since each such term is at least $-2s^{k-1} k^2/n^{k}$, their weighted average is also at least this large.
By taking expectation over $S$, we find that
\begin{equation*}
{\E[ \delta^{m1}X| F\in \SAT] }
\geq - \frac{2k^2}{n^k}\cdot\E[S^{k-1}|F\in \SAT].
\end{equation*}
Our aim is to bound $P_{m}'(b)$ in terms of $M:=\E[S^k|F\in \SAT]$. By Jensen's inequality and the convexity of $z\mapsto z^{\frac{k}{k-1}}$,
\[\E[S^{k-1}|F\in \SAT]\leq  (\E[S^{k}|F\in \SAT])^\frac{k-1}{k}=M^{\frac{k-1}{k}},\]
and plugging that into the previous inequality leads to 
\[
\E[ \delta^{m1}X| F\in \SAT] 
\geq -\frac{2k^2}{n^k}M^{\frac{k-1}{k}}.
\]
Recalling that $\frac{dP_m}{dp} = mk\E[ \delta^{m1}X]$ and $\E[ \delta^{m1}X] = P_{m-1}\E[ \delta^{m1}X|F\in \SAT]$,
\begin{equation*}
\frac{dP_m}{dp}
\geq -P_{m-1} \frac{2mk^3}{n^k}M^{\frac{k-1}{k}}.
\end{equation*}
Recalling also that $\frac{dP_m}{dp}=-\frac{dP_m}{db}$ (since $p=\frac{1}{2}-b$) gives the desired result.
\end{proof}

\noindent \begin{proof}[Proof of \cref{ineq:Pdiff-to-M}]
The difference $P_{m-1}-P_{m}$ is the probability  of the event $A:=\{\Phi_m\notin \SAT, \Phi_{m-1}\in \SAT\}$.
Since we want to lower bound this, let $B:=\{C_m \textrm{ has same sign on all variables}\}$ and consider the probability of $A\cap B$.

Conditional on $\Phi_{m-1}$, by \cref{spinevar:sat} the event $B$ occurs iff the $k$ variables in $C_m$ are either all negated and among the $s_+$ positive spine variables, or all pure and among the $s_-$ negative spine variables.
The probability of this is
\[
\frac{\binom{s_+}{k}+\binom{s_-}{k}}{\binom{n}{k}}\cdot \frac{p^k+(1-p)^k}{2}
\]
Since the function $p\mapsto p^k$ is convex on $[0,1]$, the second factor is at least $2^{-k}$.
Note also that $t\mapsto \binom{t}{k}$ is a convex function on the non-negative integers, 
%Note also that since $t\mapsto \binom{t}{k}$ is an increasing function, $\binom{s_+}{k}+\binom{s_-}{k}\geq \max\binom{s_\pm}{k}\geq \binom{s/2}{k}$, 
%where $s:=s_++s_-$. Hence the first factor is at least ${f(s):=\binom{s/2}{k}/\binom{n}{k}}$.
%And since the function $t\mapsto \binom{t}{k}$ is convex for $t\geq k$,
and hence the first factor is at least ${f(s):=2\binom{s/2}{k}/\binom{n}{k}}$,
where $s:=s_++s_-$.
The standard bounds ${(t/ek)^k\leq\binom{t}{k}\leq (t/k)^k}$ are valid for $t\geq k$, so for ${s\geq 2k}$ we have that $f(s)\geq 2(s/2en)^{k}$.
If instead $s<2k$, $f(s)\geq 0 > \lambda\cdot (s^k-(2k)^k)$ for any $\lambda>0$. In either case,
\[
{f(s)\,\geq\,   2(2en)^{-k}\cdot \big(s^k-(2k)^{k}\big)}.
\]
Then, by taking expectation over satisfiable $\Phi_{m-1}$ (and hence over $S$),
\begin{align}
\Pr(A\cap B|\Phi_{m-1}\in \SAT)
&\geq \E[2^{-k}f(S) |\Phi_{m-1}\in \SAT]
\nonumber
\\
&\geq
2(4en)^{-k}\cdot  \left(\E[S^k |\Phi_{m-1}\in \SAT]-(2k)^k\right) 
\nonumber
\end{align}
Noting that the left-hand side is at most $\Pr(A|\Phi_{m-1}\in \SAT) = (P_{m-1}-P_m)/P_{m-1}$ and solving for $M=\E[S^k\big|\Phi_{m-1}\in \SAT]$) gives the desired inequality.
\end{proof}

\noindent \begin{proof}[Proof of \cref{prop:roughbounds}]
The first inequality is an easy corollary of \cref{lemma:2SAT-lb}: Given $F:=\Phi_m(n,k,p)$, $k\geq 3$, uniformly at random remove all but $2$ literals from each clause to get a $2$-SAT formula $\tilde F$. Any satisfying assignment to $\tilde F$ is also a satisfying assignment to $F$, and $\tilde F$ is a $p$-polarized $2$-SAT formula. Hence $\Pr(\Phi_k\in \SAT) \geq \Pr(\Phi_2\in \SAT)$, and by \cref{lemma:2SAT-lb} this probability is at least ${1-30/(n\eps^3)}$ for ${m<(1-\eps) n}$.
For the second inequality we will upper bound the expected number of satisfying assignments.

\begin{claim}
\label{claim:qi-bound}
Let $q_i:=\Pr(C(\sigma)=-1)$ be the probability that a $\sigma\in \{\pm 1\}^n$ with $i$ coordinates set to 'TRUE' does not satisfy a random $p$-polarized $k$-clause $C$. Then $q_i \geq  2^{-k}-k^2/n$
\end{claim}
\vspace{-1em}
\noindent \begin{proof}[Proof of claim]
Let $K,\tilde K$ be random $k$-tuples $(v_1,\ldots v_k)$ and $(\tilde v_1,\ldots \tilde v_k), $ chosen uniformly at random from the set $\{1\,\ldots ,n\}$, \emph{without} and \emph{with} replacement respectively.
By \cref{claim:with-or-without-replacement}, we can couple $K,\tilde K$ such that they are equal with probability at least $1-k^2/n$. Like in \cref{setup}, let the random variable $B$ equal $\pm 1$ with probability $\frac{1}{2}$, and $P_j=1$ with probability $p$ (and $-1$ otherwise). Let $z_j := x_{v_j}$ and $z_j:= x_{\tilde v_j}$ if $BP_j=1$, but $z_j := \neg x_{v_j}$ and $\tilde z_j := \neg x_{\tilde v_j}$ otherwise. Finally, let $C:=z_1\vee\ldots \vee z_k$ and $\tilde C:=\tilde z_1\vee\ldots \vee \tilde z_k$ be $k$-clauses.  

Conditional on $B=1$, the events $\{\tilde z_j(\sigma)=-1\}_{j=1}^k$ are independent and each occur with the same probability $\rho(\sigma)$. Conditional on $B=-1$, they also occur independently, but with probability $1-\rho$.
Hence the probability that $\sigma$ does not satisfy $\tilde C$ is $\frac{1}{2}(\rho^k+(1-\rho)^k) \geq 2^{-k}$.
But $C=\tilde C$ with probability at least $1-k^2/n$, so the probability that $\sigma$ does not satisfy $C$ is at least $2^{-k}-k^2/n$.
\end{proof}

Let $m> -n/\log_2(1-2^{-k})+t$. The logarithm of the expected number of satisfying assignments to $\Phi_m(n,k,p)$ is then equal to
\begin{align*}
&\log_2\left(\sum_{i=1}^n\binom{n}{i}(1-q_i)^{m}\right)
\leq \log_2\left(2^n \cdot (1-2^{-k}+k^2/n)^{m}\right)
\\
=&n \cdot \Big( \underbrace{1-\frac{\log_2 (1\!-\!2^{-k}\!+\!k^2/n)}{\log_2(1-2^{-k})}}_{(i)}\Big)+ \underbrace{t\log_2(1-2^{-k}+k^2/n)}_{(ii)}.
\end{align*}
By doing a first order Taylor expansion of the increasing and convex function ${x\mapsto 1-\log_2(1-2^{-k}+x)/\log_2(1-2^{-k})}$ at $x=0$, we see that $(i)\leq 2^k k^2/n$ for large $n$. The first term above is therefore at most ${2^k k^2}$. For $n>2^{k+1} k^2$, the second term $(ii)$ is at most ${t\log_2(1-2^{-k-1})<-2^{-k-2} t}$. The second inequality of the proposition follows.
\end{proof}

\section{Proof of \cref{thm:2SAT}}
\label{section:k2}
For the classical $2$-SAT problem, the threshold value of $c = 1$ was established by Chvátal \& Reed~\cite{chvatal-reed} by exploiting some of the structure specific to $2$-SAT. 
Our proof is fairly similar to theirs, but with some minor complications. 

A $2$-clause is of the form $x\vee y$, which is logically equivalent the implication $\neg x \Rightarrow y$, and also to the implication $\neg y \Rightarrow x$. Thus, a $2$-SAT formula with $m$ clauses on $n$ variables can be represented by a digraph $G$ with $2m$ arcs on the following $2n$ vertices: $\{x_1,\ldots, x_n, \neg x_1,\ldots, \neg x_n \}$. We call any directed cycle in G a \emph{bicycle}\footnote{Our notation follows that of~\cite{chvatal-reed} and later papers on 2-SAT.}  if it contains both $x_i$ and $\neg x_i$ for some $i$. It is clear that the formula is unsatisfiable if the digraph contains a bicycle, because from a bicycle the contradiction $x_i\Leftrightarrow \neg x_i$ can be derived. Less obviously, this is an 'if and only if'-condition.
%\vspace{-1em}
\begin{lemma}[Aspvall, Plass \& Tarjan~\cite{2SAT:iffconditionbicycle}]
A $2$-SAT instance is satisfiable if and only if there is no bicycle in the associated digraph of implications.
\end{lemma}
\noindent We will establish upper and lower bounds on the satisfiability threshold by applying the second and first moment method to certain structures related to bicycles.  
\begin{define}
A \emph{unicycle} $($short for unique bicycle $)$ is an even length bicycle with a unique repeated variable $x$ and with that variable occurring precisely twice, once as $x$ and once as $\neg x$, at diametrically opposed points along the cycle.

A \emph{pretzel} is a directed path $\ell\to \ell_1\to \ldots \to \ell_t\to\ell'$, where each $\ell_i$ is a literal on $z_i$, the $z_i$'s are distinct, and $\ell,\ell'$ are literals on some variables $z,z'\in \{z_1, \ldots ,z_t\}$.

A pair of edges is \emph{conjugate} if the corresponding implications are contrapositives of one another (i.e. both implications arise from the same clause). 
\end{define}
\begin{lemma} \label{lemma:cycleequiv}
For directed graphs $G$ we have that:
\[(\exists \textrm{ unicycle} \in G) \Rightarrow (\exists \textrm{ bicycle} \in G) \Rightarrow (\exists \textrm{ pretzel} \in G)\]
\end{lemma}
\noindent\begin{proof}
The first implication is clear: any unicycle is a bicycle.
For the second one, let $C\subseteq G$ be a bicycle.
Pick a maximal path $P\subset C$ such that the literals in $P$ are on disjoint variables. We can then write $P$ as $\ell_1 \to \ell_2 \to  \ldots \to \ell_t$ for some $t>0$, where each $\ell_i$ is a literal on a variable $z_i$, and the $z_i$'s are distinct.

Let $\ell,\ell'$ be the unique literals in $C$ such that $\ell \to \ell_1$ and $\ell_t\to \ell'$, or in other words the immediate predecessor and successor to $P$. It follows from the maximality of $P$ that $\ell,\ell'$ are literals on some variables in the set $\{z_1,\ldots ,z_t\}$. (If $\ell$ or $\ell'$ were a literal on a variable not in this set, we could extend $P$ to a larger path of literals on disjoint variables.) Hence the path $\ell\to \ell_1 \to \ell_2 \to  \ldots \to \ell_t \to \ell'$ is a pretzel.
\end{proof}

We will use the first moment method to show that for any $0\leq p \leq 1$ there is no pretzel (w.h.p.) when $m < (1-\eps)n$ (and hence there is no bicycle).
We will use the second moment method to show that for $p=0$ there is a unicycle (w.h.p.) when $m > (1+\eps)n$ (and hence there is a bicycle). Then \cref{prop:increasingthreshold} will imply that the same is true for every $0\leq p \leq 1$.
\begin{prop}
\label{lemma:2SAT-lb}
For any $\eps=\eps(n)$, if $m<(1-\eps)n$ then the probability that there exists a pretzel is at most $30/(n\eps^3)$.
\end{prop}
\noindent\begin{proof}
We will begin by upper bounding the expected number of pretzels with $t+1$ edges.
Note that no arc in a pretzel is conjugate to another arc in it,  so they all correspond to distinct $2$-clauses. Furthermore, a pair of variables $x_i,x_j$ can occur in at most one $2$-clause (since we pick such pairs \emph{without} replacement).% when we generate the random $2$-SAT formula).

There are at most $n^t$ ways to choose the $t$ variables $z_i$, and then at most $t^2$ ways to choose $z$ and $z'$ from the $t$-set $\{z_1,\ldots z_t\}$, so there are at most $n^t t^2$ ways to choose the variables in a pretzel.

The probability that each of the $t+1$ pairs of variables $zz_1, z_1 z_2, \ldots, z_{t-1}z_t, z_t z'$ occurs in some $2$-clause in the $2$-SAT formula is at most $(m/\binom{n}{2})^{t+1}$.
Conditioned on these $t+1$ pairs occurring as $2$-clauses, what is the probability that the signs of variables in all of them are such that the corresponding arcs in the digraph form a directed path? In other words, that the $2$-clauses form implications ${\ell\Rightarrow \ell_1\Rightarrow \ldots \Rightarrow \ell_t\Rightarrow \ell'}$, where the $\ell_i$'s are literals on the variables $z_i$.

We reveal the signs of variables in clauses in order:
First, condition on the event that $z_1$ occurs pure in the first clause (whose variables are $z$ and $z_1$), i.e. $\ell_1=z_1$ so that the clause is $\neg\ell\vee z_1$, for either $\ell=z$ or $\ell=\neg z$. This is equivalent to the implication $\ell \Rightarrow z_1$.

Then, in order for the next clause (with variables $z_1, z_2$) to form an implication of the form $z_1 \Rightarrow \bullet$, the variable $z_1$ must occur negated in it, but the sign of $z_2$ is not constrained. In other words, the second clause must be either $\neg z_1 \vee z_2$ or $\neg z_1 \vee \neg z_2$. The probability of this event is $\frac{1}{2}$.
But if we instead condition the first clause being $\ell \Rightarrow \neg z_1$, the probability that the next clause is $\neg z_1\Rightarrow \bullet$ is also $\frac{1}{2}$.

Similarly, for each subsequent clause (with variables $z_i, z_{i+1}$), the sign of $z_i$ in it must be the opposite of the sign that $z_i$ had in the previous clause, which happens with probability $\frac{1}{2}$ independently of previous clauses. The probability that all these $t$ events occur is $\left(\frac{1}{2}\right)^{t}$. 
The expected number of pretzels is thus at most
\begin{align*}
\E[\#\textrm{ pretzels}]\leq\, &\sum_{t=1}^{2n} n^{t} t^2 \cdot \left( \frac{m}{\binom{n}{2}}\right)^{t+1}\left(\frac{1}{2}\right)^{t}
\\
=\,
&\frac{2}{n}\cdot\sum_{t=1}^{2n} t^2  \left(\frac{m}{n-1} \right)^{t+1}
\\
<\,
&\frac{2}{n}\left(\frac{n}{n-1}\right)^{2n}\cdot\sum_{t=1}^{2n} t^2  \Bigg(\frac{m}{n} \Bigg)^{t+1}.
\end{align*}
By the assumption on $m$, $\frac{m}{n}< 1-\eps$. But then the sum above is at most
 \[\sum_{t=1}^{\infty} t^2 (1-\eps)^{t+1} = \frac{(1-\eps)^2(2-\eps)}{\eps^3}< \frac{2}{\eps^3},\]
so that the expected number of pretzels is at most
\[
\frac{2}{n}\left(\frac{n}{n-1}\right)^{2n}\cdot \frac{2}{\eps^3} < 4e^2 \cdot \frac{1}{n\eps^3}.
\]
Hence there exists a pretzel when ${m<(1-\eps)n}$ with probability $\leq 30/(n \eps^3)$.
\end{proof}

\begin{prop}\label{lemma:2SAT-ub}
Let $p=0$ and $\frac{1}{2}>\eps \gg n^{-0.1}(\ln n)^{0.9}$. Then if ${m>(1+\eps)n}$ there exists a unicycle (with high probability).
\end{prop}
\noindent\begin{proof}
By assumption, $(n \ln n)^{0.9}\ll \eps n$. Let $t=t(n)$ be a sequence of odd integers such that $(n \ln n)^{0.9} \ll  t^9\ll \eps n$.
We will consider only unicycles on precisely $2t$ vertices and show that at least one such unicycle is present with high probability.

The unicycles are much more constrained for the fully polarized $2$-SAT problem than for the classical $2$-SAT problem. Since clauses of the form $x_i \vee \neg x_j$ occur with probability $0$,  there will be no arcs of the form $x_i \Rightarrow x_j$ or $\neg x_i \Rightarrow \neg x_j$ in the digraph of implications. In other words, the digraph becomes bipartite, with the vertices partitioned into the sets $V_+:=\{x_1,x_2,\ldots x_n\}$ and $V_-:=\{\neg x_1,\neg x_2,\ldots \neg x_n\}$.  The proof of Theorem 4 of~\cite{chvatal-reed} uses the second moment method applied to the number of unicycles on the complete digraph, and we will describe how to adapt the proof to instead count the number of unicycles on the complete \emph{bipartite} digraph.

First, note that the total number of possible clauses in the fully polarized model is $2\binom{n}{2}$, compared to $4\binom{n}{2}$ in the original model. Note also that a path from $x_i$ to $\neg x_i$ has an odd number of edges, so it is necessary that $t$ is odd.

Let $(n)_t$ denote the falling factorial $n(n-1)\ldots(n-t+1)$. A directed path $\ell_1   \Rightarrow  \ldots \Rightarrow \ell_t $, where each $\ell_i$ is a literal on the variable $z_i$, is determined by the underlying sequence of variables $z_i$ together with whether $\ell_1=z_1$ or $\ell_1=\neg z_1$, so there are $2\cdot (n)_{t}$ paths of length $t$ in the complete bipartite digraph. Note that $1> (n)_{t}/n^{t}> (1-t/n)^t$, and since $t\ll \sqrt{n}$ by assumption, $2\cdot (n)_{t}=(1+o(1))\cdot 2n^t$.
Similarly, the complete digraph on $2n$ vertices has $(2n)_t=(1+o(1)) \cdot (2n)^t$ paths. 

Keeping these difference in mind, the estimates in Theorem 4 of~\cite{chvatal-reed}, pages 625--626, of the first two moments of the number of unicycles on $2t$ vertices carry through with minimal alterations: Replace all powers of $2n$ with powers of $n$, every power of $2$ with a power of $1$, except in estimate $(iii)$, and replace $4\binom{n}{2}$ with $2\binom{n}{2}$ throughout the proof of said theorem.
In the end, all of these alterations cancel out, and we reach the same conclusion that as long both
\begin{align}
\label{CR:estimate1} t^9/\eps n&=o(1), \textrm{ and }  \\ \label{CR:estimate2} tn (1-\eps)^t /\eps &=o(1),
\end{align}
then there exists at least one unicycle (w.h.p.) whenever $m>(1+\eps)n$. Condition (\ref{CR:estimate1}) holds by our choice of $t$. To show that (\ref{CR:estimate2}) holds, we deal with ${tn}/{\eps}$ and $(1-\eps)^t$ separately. Note that $\eps t \gg t^{10}/n \gg \ln n$, so 
\[
(1-\eps)^t   \leq  \exp(-\eps t) =   \exp(-\omega(\ln n))=n^{-\omega(1)}.
\]
On the other hand, $tn/\eps < t^9n/\eps =o(n^2)$ by our choice of $t$. It follows that ${tn (1-\eps)^t /\eps =o(1)}$. Hence both (\ref{CR:estimate1}) and (\ref{CR:estimate2}) hold, and the theorem follows.
\end{proof}

\noindent \begin{proof}[Proof of \cref{thm:2SAT}]
\Cref{lemma:2SAT-lb,lemma:cycleequiv} together show that $\Phi_m(n,2,p)$ is satisfiable w.h.p. for any $m\leq n - \omega(n^{2/3})$.
For the upper bound, by \cref{mainthm}:
\[\Pr(\Phi_m(n,2,p) \notin \SAT) \geq \Pr(\Phi_m(n,2,0) \notin \SAT)-o(1).\]
\cref{lemma:2SAT-ub,lemma:cycleequiv} together show that $\Phi_m(n,2,0)$ is unsatisfiable w.h.p. for any $m\geq n+\omega((n \ln n)^{9/10})$, so the right hand side above is $1-o(1)$. In other words, $\Phi_m(n,2,p)$ is unsatisfiable w.h.p.\footnote{Note that  \cref{prop:roughbounds} depends on the lower bound on the $2$-SAT threshold (\cref{lemma:2SAT-lb}), while the upper bound on the $2$-SAT threshold (\cref{lemma:2SAT-ub}) depends on \cref{mainthm}, which in turn depends on \cref{prop:roughbounds}. But there is no circularity of reasoning here, since the upper and lower bounds on the $2$-SAT threshold do not depend on each other.}
\end{proof}

\bibliographystyle{abbrv}
\bibliography{satpapers}
\end{document}